\documentclass[review]{elsarticle}

\usepackage{lineno,hyperref,lscape}

\modulolinenumbers[5]
\usepackage{amsmath,amssymb}
\journal{arXiv}

\begin{document}

\begin{frontmatter}

\title{Fractional Calculus Approach to Logistic Equation and its Application}

\author{Jignesh P. Chauhan$^{1}$, Ranjan K. Jana$^{2}$, Pratik V. Shah$^{3}$, Ajay K. Shukla$^{4}$}

\address{$^{1,2,4}$Department of Applied Mathematics \& Humanities, S. V. National Institute of Technology, Surat- 395007\\
Email: impulse.nit07@gmail.com$^{1}$, rkjana2003@yahoo.com$^{2}$, ajayshukla2@rediffmail.com$^{4}$}
\address{ $^{3}$Department of Applied Sciences, C. K. Pithawalla College of Engineering and Technology,

Surat-395007, India. \\
Email: pratikshah8284@yahoo.co.in$^{3}$}

\begin{abstract}
In this paper, we propose a solution of fractional logistic equation by using properties of Mittag-Leffler function.
\end{abstract}

\begin{keyword}
Mittag-Leffler function;  Caputo fractional derivative; Logistic equation.\\
\textit{2010 AMS Classification}: 34A08, 26A33, 33E12
\end{keyword}
\end{frontmatter}

\linenumbers

\section{Introduction}
\noindent In recent years applications of fractional calculus have been investigated extensively. Considerable amount of work has been done in area of fractional differential equations and many analytical and numerical methods were developed and employed for obtaining the solution reported by Mathai \textit{et al.} \cite{mathai2008special}. It has been found that in some cases fractional calculus is more accurate than classical calculus to describe dynamic behavior of real world physical systems.\\
Exponential function which arises by solving differential equation plays an important role for describing growth and decay in many physical applications. In fractional ordered differential equation, exponential function loses its properties to describe the solution and Mittag-Leffler function is used as its substitute. \\
A non linear differential eqauion of population growth model was first published by Verhulst \cite{verhulst1838notice}, subsequently known as logistic equation,
\begin{align}\label{le}
{{du\left( t \right)} \over {dt}} = k~u(t) {\left( {1 - u\left( t \right)} \right)} , ~~~ t \ge 0,
\end{align}
whose exact closed form is given by,
\begin{align}\label{les}
u\left( t \right) = {{{u_0}} \over {{u_0} + \left( {1 - {u_0}} \right)\exp ( - kt)}},
\end{align}
where $u_0$ is the initial state when time $t=0$.
This equation often arises while modeling ecology, neural networks, epidemics, Fermi distribution, economics, sociology etc. So, we are motivated to study fractional logistic equation by generalizing \eqref{le} to its non-integer part.

In 2015, Bruce J. West \cite{west2015exact} considered fractional form of non-linear logistic equation  in following form,
\begin{align}\label{bfle}
D_t^\alpha \left[ {u\left( t \right)} \right] = {k^\alpha }u\left( t \right)\left[ {1 - u\left( t \right)} \right].
\end{align} 
He found solution using Carleman embedding technique as,
\begin{align}\label{bfles}
u\left( t \right) = \sum\limits_{n = 0}^\infty  {\left( {{{{u_0} - 1} \over {{u_0}}}} \right)^n {E_\alpha }\left( { - n{k^\alpha }{t^\alpha }} \right)} ,~~~ t \geq 0.
\end{align}
Whereas in 2016, Area \textit{et al.} \cite{area2016note}, has shown that the solution given in \eqref{bfle} is not an exact solution of fractional logistic equation. Further, Ortigueira \textit{et al.} \cite{ortig} expressed exact solution of fractional logistic equation in terms of fractional Taylor series.

In this paper, we propose a new approach in light of Jumarie \cite{jumarie2009table} concept to obtain solution of fractional logistic equation.

 \section{Definitions}
\subsection{Mittag-Leffler Function} 
The generalized form of exponential function, Mittag-Leffler function \cite{gorenflo2014mittag}, is named after Swedish mathematician G{\"o}sta Mittag-Leffler is given by,\begin{align}
{E_\alpha }\left( z \right) = \sum\limits_{k = 0}^\infty  {{{{z^k}} \over {\Gamma \left( {\alpha k + 1} \right)}}} ,~~~~\alpha  > 0.
\end{align}
 \subsection{Mittag-Leffler function in two parameters}
\noindent The \textit{Mittag-Leffler function} \cite{gorenflo2014mittag} with two parameters is defined as,  
\begin{equation} \label{1.2} 
E_{\alpha ,\beta }\left(z\right)=\ \sum^{\infty }_{n=0}{\frac{z^n}{\Gamma \left(\alpha n+\beta \right)}\ \ }, 
\end{equation} 
where $\alpha ,\beta \in C;Re\left(\alpha \right)>0\ ,\ Re(\beta )>0. $\\
Due to direct involvement in generalization of ordinary differential equation to its non-integer order, Mittag-Leffler function is found very useful in many areas of science and engineering. 
\subsection{Caputo's Fractional Derivative}
\textit{Caputo's} \cite{podlubny1998fractional} definition of fractional derivative is given by 
\begin{align}\label{11.1.2}
 \,_0^CD_t^\alpha f\left( t \right) = {1 \over {\Gamma \left( {n - \alpha } \right)}}\int_0^t {{{{f^n}\left( \tau  \right)} \over {{{\left( {t - \tau } \right)}^{\alpha  - n + 1}}}}d\tau ,} 
\end{align}
where $\alpha \in R$ is order of fractional derivative, $n-1<\alpha \leq n$ and $n \in N= \{ 1,2,3,\dots \}$, ${f^n}\left( \tau  \right) = {{{d^n}} \over {d{t^n}}}f\left( \tau  \right)$ and $\Gamma (.)$ is Euler Gamma function.\\[.2cm]

\section{Main Results}
\textbf{Logarithm Function with Mittag-Leffler Function ($E_\alpha )$ as the base}\\
Let $E_\alpha : \mathcal{R} \rightarrow \mathcal{R^+}$ (one-one\footnote{\textbf{This is important to note that $E_\alpha$ is not one to one in general and hence the $\log$ function of $E_\alpha$ is not always defined.}} and onto , if exists) has inverse function $L_\alpha (\log _{E_\alpha})$, we define $L_\alpha$ as
\begin{align} \label{c1}
\left.
\begin{matrix}
{E_\alpha }\left[ {{L_\alpha }\left( y \right)} \right] = y,~& where &y > 0 \\
 {L_\alpha }\left[ {{E_\alpha }\left( x \right)} \right] = x, &where &x \in \mathcal{R}\\
{E_\alpha }\left( x \right) = y &or &{L_\alpha }\left( y \right) = x \\
{E_\alpha }\left( x \right) \to  + \infty &as &x \to  + \infty   \cr 
   {E_\alpha }\left( x \right) \to 0\,\hspace*{.3cm}&as &x \to  - \infty   \cr 
   {L_\alpha }\left( x \right) \to  + \infty &as &x \to  + \infty   \cr 
   {L_\alpha }\left( x \right) \to  - \infty  &as &x  \to 0 \hspace*{.5cm} \cr
\end{matrix} \right\rbrace,  ~~~where ~~Re(\alpha)>0.
\end{align}
Mittag-Leffler function $E_\alpha (x)$ is differentiable and its inverse function $L_\alpha (\log _{E_\alpha})$ is also differentiable, this is important to note that  $E'_\alpha (x)\neq E_\alpha (x)$.
However many authors have proved some properties of Mittag-Leffler function by considering ${E_\alpha }\left( {a{{\left( {{x_1} + {x_2}} \right)}^a}} \right) = {E_\alpha }\left( {a{x_1^a}} \right){E_\alpha }\left( {a{x_2^a}} \right)$, where ${x_{1,}}~~{x_2} \ge 0$, $a$ is real constant and $\alpha >0$. Neito \cite{nieto2010maximum} gave many interesting results based on  Mittag-Leffler function. Jumarie \cite{jumarie2009table} also gave various definitions of fractional ordered differentiation and integration. Many important results have also been  reported by Gorenflo \textit{et al.} \cite{gorenflo2014mittag} in their book.  \\
We are motivated to establish some new results on Mittag-Leffler function in the form of following preposition, \\

\noindent \textbf{Proposition:} \textit{Let $u=E_\alpha (x_1)$ and $v=E_\alpha (x_2)$ , where $\alpha >0 $ with above conditions \eqref{c1}. Then , }
\begin{enumerate}
\item[(i)] $E_\alpha (x_1 \oplus x_2)=E_\alpha (x_1) \odot E_\alpha (x_2)$
\item[(ii)] $\log_{E_\alpha} \left( {u \odot v}\right)=x_1 \oplus x_2=\log_{E_\alpha} \left( {u}\right) \oplus \log_{E_\alpha} \left( {v}\right)$ 
\item[(iii)] $\log_{E_\alpha} \left( {u ~{\tiny \raisebox{.1pt}{\textcircled{\raisebox{-.5pt} {$\div$}}}}~ v}\right)=x_1 \ominus x_2=\log_{E_\alpha} \left( {u}\right) \ominus \log_{E_\alpha} \left( {v}\right)$ 
\end{enumerate}
\textbf{Proof:}
\begin{enumerate}
\item[(i)] From \eqref{c1}, we consider 
\begin{align} 
& x_1 \oplus x_2 = L_\alpha (u \odot v), 
\end{align}
then
\begin{align}
 & E_\alpha \left( x_1 \oplus x_2 \right) = u \odot v = E_\alpha \left( x_1 \right) \odot E_\alpha \left( x_1 \right) \label{l},
\end{align}
when $\alpha \rightarrow 1$, then \eqref{l} reduces to 
\begin{align}
E_1 \left( x_1 \oplus x_2 \right)= E_1 \left( x_1 \right) \odot E_1 \left( x_1 \right). 
& \nonumber
\end{align}
This is same as,
\begin{align}
\exp \left( x_1 + x_2 \right)&= \exp \left( x_1 \right) \cdot \exp \left( x_1 \right) .
\end{align}
\item[(ii)] Now we consider 
\begin{align}
u\odot v=E_\alpha \left( x_1 \oplus x_2 \right),
\end{align}
 taking logarithm (base $E_\alpha $) on both of the sides,
\begin{align}
\log_{E_\alpha} \left( u\odot v \right)&=x_1 \oplus x_2 \nonumber \\ &= 
\log_{E_\alpha} (u) \oplus \log_{E_\alpha} (v). \label{l2}
\end{align}
On taking limit $\alpha \rightarrow 1 $ then \eqref{l2} reduces to 
\begin{align}
\log_{E_1} \left( u\odot v \right)= 
\log_{E_1} (u) \oplus \log_{E_1} (v).
\end{align}
Which is same as,
\begin{align}
\log_e (u\cdot v)=\log_e (u)+\log_e (v).
\end{align}
\item[(iii)]  Similarly, we can easily show 
\begin{align}
\log_{E_\alpha} \left( {u ~{\footnotesize \raisebox{.12pt}{\textcircled{\raisebox{-.2pt} {$\div$}}}}~ v}\right)=x_1 \ominus x_2=\log_{E_\alpha} \left( {u}\right) \ominus \log_{E_\alpha} \left( {v}\right). \label{l3}
\end{align}
On taking $ limit ~~ \alpha \rightarrow 1$ then \eqref{l3} reduces to
\begin{align}\label{odiv2}
\log_{E_1} \left( {u ~{\footnotesize \raisebox{.12pt}{\textcircled{\raisebox{-.2pt} {$\div$}}}}~ v}\right)=x_1 \ominus x_2=\log_{E_1} \left( {u}\right) \ominus \log_{E_1} \left( {v}\right) ,
\end{align}
or,
\begin{align}
\log_e \left( \frac{u}{v} \right)=\log_e (u)- \log _e(v).
\end{align}
\end{enumerate}

\textit{Note:} It is interesting to see that operators $\oplus,~\ominus,~\odot ~and ~{\footnotesize \raisebox{.12pt}{\textcircled{\raisebox{-.2pt} {$\div$}}}}~$  behave very closed to  traditional addition, subtraction, multiplication and division operators for different values of $\alpha, ~x_1~and~x_2$. One can easily verify above proposition numerically by substituting differnt values for $\alpha$, $0< \alpha \leq 1$, from Table.\ref{table1} and Table.\ref{table2}, where $\log_{E_\alpha}$ is logarithm function with Mittag-Leffler function as base.
Some properties of logarithmic function with base Mittag-Leffler function $E_\alpha(=E_\alpha (1))$ seems to be similar as logarithmic function with base $e$. The graph of $\log_{E_\alpha}$ for different choices of $\alpha $ is shown in Figure.\ref{Fig1}.\\

\subsection{Solution of fractional logistic equation}

\noindent We constitute \eqref{le} in fractional differential equation as,  
\begin{align}\label{3.4}
& \frac{d^\alpha u}{dt^\alpha}= k^\alpha u\left( {1 - u} \right),~where~ 0<\alpha \leq 1,~~u=u(t)
\end{align}
On writing \eqref{3.4} in following manner,
\begin{align} \label{3.4_2} & \int {{{d^\alpha {u}} \over u} \oplus \int {{{d^\alpha {u }} \over {1 - u}} = k^\alpha \int {d{t^\alpha }} } } .
\end{align}
Jumarie \cite{jumarie2009table} clarified that, some formulas do not hold for the classical Riemann-Liouville definition, but can be applied with the modified Riemann-Liouville definition.\\
Further using the proposition for obtaining the solution of \eqref{3.4_2}, we find, 

\begin{align}\label{3.5_1}
{{{\log }_{{E_\alpha }}}u \ominus {{\log }_{{E_\alpha }}}\left( {1 - u} \right)}= {\frac{k^\alpha}{\Gamma(2-\alpha)}} \int t^{1-\alpha}dt^{\alpha} \oplus C, 
\end{align}
where $C$ is the integration constant. On using \eqref{c1} and \eqref{odiv2}, we get
\begin{align}
 \label{3.5_2} \frac{u(t)}{1-u(t)}=CE_\alpha \left[ {\frac{k^\alpha }{\Gamma(2-\alpha)}} \int t^{1-\alpha}dt^{\alpha} \right] ,
\end{align}
where $0 < \alpha  \leq 1$ and on putting $\alpha=1$ we get ${\log _{{E_\alpha }}}x = {\log _e}x$.  Here we also clarify that the numerical value of $u(t) ~{\footnotesize \raisebox{.12pt}{\textcircled{\raisebox{-.2pt} {$\div$}}}}~ (1-u(t)) $ is approximately very close to $\frac{u(t)}{1-u(t)}$.
 
\noindent 
Initially, when time $t=0$, we write,
$u(0)=u_0$ and  $C=\frac{u_0}{1-u_0}$.
Now substituting the value of $C$ in \eqref{3.5_2}, this gives 
\begin{align}\label{3.6}
&\frac{u(t)}{1-u(t)}=\frac{u_0}{1-u_0}E_\alpha \left[ {\frac{k^\alpha}{\Gamma(2-\alpha)}} \int t^{1-\alpha}dt^{\alpha} \right]. 
\end{align}
On further simplification, we arrive at
\begin{align}\label{jle}
 & u(t)=\frac{1}{1+\frac{1-u_0}{u_0} \left[ E_\alpha \lbrace {\frac{k^\alpha}{\Gamma(2-\alpha)}} \int t^{1-\alpha}dt^{\alpha} \rbrace \right]^{-1}},
\end{align}
On setting $\alpha = 1$, equation \eqref{jle} reduces to \eqref{les}.

\newpage

\section{Application}
\textbf{Epidemic Model:}
Epidemiology is concerned with the spread of disease and its effect on people. This in itself encompasses a range of disciplines, from biology to sociology and philosophy, all of which
are utilized to a better understanding and containing of the spread of infection. Classical Epidemics transmission models have been used to interpret the spread of epidemic without immunity after an infective period, such as encephalitis, gonorrhoea, etc.  
The earliest account of mathematical modeling for spread of disease was carried out in 1766 by Daniel Bernoulli. It can be used to explain the change in the number of people needing medical attention during an epidemic.
One of the work related to mathematical theory of epidemics was given by Kermack and McKendrick \cite{kermack1927contribution}.
Based on some mathematical assumptions, epidemics can be modeled mathematically in order to study the severity and prevention mechanism by which diseases spread. This also helps to predict the future course of an outbreak and to evaluate strategies to control an epidemic.\\
Here we assume that the whole population $N$ is divided into three sections, i.e. $S$, the number of susceptible; $I$, the number of infected and $R$, the number of recovered during an epidemic. This model assumes that the total population remains the same with closed demography, i.e. there is no birth and no natural death. Any disease related to death, however, can be included in $R$. Recently many authors have studied and developed epidemic models by using methods of bifurcation \cite{kermack1927contribution}. The ultimate goal is to model the issue of saturated susceptible population, the time delay of infected to become infectious, the stability of equilibrium solutions. \\
As the first step in the modeling process, we identify the independent and dependent variables. The independent variable is time  $t$,  measured in days. We consider two related sets of dependent variables.\\
$S = S(t)$	is the number of susceptible individuals,\\
$I = I(t)$	is the number of infected individuals, and\\
$R = R(t)$	is the number of recovered individuals.\\
Considering total population during epidemic spread as $N$, we have, 
\begin{align} \label{1}
S\left( t \right) + I\left( t \right) + R\left( t \right) = N\left( t \right)
\end{align} 
Here we consider three basic epidemic models with some reasonable assumptions, the epidemic models are:
\subsection{\textbf{SI Model}}
In this model, we assume that total population consist of only susceptible and infected people.  Therefore, total population is given by \cite{kapur1988mathematical},$$ S(t)+I(t)=N(t).$$
The epidemic starts with $I_0$ (at time $t=0$) infected persons. The rate at which susceptible changes to infected with respect to time $t$ is given by  
\begin{align} \label{3.1}
{{dS} \over {dt}} =  - \beta I\left( t \right)S\left( t \right), ~~~t > 0,
\end{align}
where $\beta$ is the positive integer and initially at time $t=0,~I(0)=I_0$. 
Also, the rate of change of infected is given by 
\begin{align}\label{3.2}
{{dI} \over {dt}} &= \beta I\left( t \right)\left( {N - I} \right), ~~~ I\left( 0 \right) = {I_0} ~~ and  \nonumber \\  
 {{dS} \over {dt}} &=  - \beta S\left( t \right)\left( {N - S} \right), ~~~ S\left( 0 \right) = N - {I_0} .
\end{align}
\subsubsection{\textbf{Fractional Differential Equation  for SI Model:}}
Here we are motivated to study the fractional differential equation model for epidemic and obtain the solution for the said model. \\
On writing \eqref{3.2} with arbitrary order $\alpha$, as,  
\begin{align}
\label{3.4}& \frac{d^\alpha I}{dt^\alpha}= \beta I\left( {N - I} \right),~where~ 0<\alpha \leq 1\\  \label{3.4_2} \Rightarrow & \int {{{d^\alpha {I}} \over I} + \int {{{d^\alpha {I }} \over {N - I}} = N\beta \int {d{t^\alpha }} } } .
\end{align}
On using the \eqref{odiv2}, the solution of \eqref{3.4_2} is given by, 
\begin{align}\label{3.5}
{{{\log }_{{E_\alpha }}}I \ominus {{\log }_{{E_\alpha }}}\left( {N - I} \right)}= {\frac{N\beta}{\Gamma(2-\alpha)}} \int t^{1-\alpha}dt^{\alpha} + c, \\
\Rightarrow \frac{I}{N-I}=CE_\alpha \left[ {\frac{N\beta}{\Gamma(2-\alpha)}} \int t^{1-\alpha}dt^{\alpha} \right] ,
\end{align}
where $0 < \alpha  < 1$ and $\mathop {\lim }\limits_{\alpha  \to 1} \,\,{\log _{{E_\alpha }}}x = {\log _e}x.$
 
\noindent 
Initially, when time $t=0$ we have,
$I(0)=I_0$, we get $C=\frac{I_0}{N-I_0}$, i.e. 
\begin{align}\label{3.6}
&\frac{I}{N-I}=\frac{I_0}{N-I_0}E_\alpha \left[ {\frac{N\beta}{\Gamma(2-\alpha)}} \int t^{1-\alpha}dt^{\alpha} \right] \\ 
\Rightarrow & I=\frac{N}{1+\frac{N-I_0}{I_0} \left[ E_\alpha \lbrace {\frac{N\beta}{\Gamma(2-\alpha)}} \int t^{1-\alpha}dt^{\alpha} \rbrace \right]^{-1}}
\end{align}
when $\alpha \rightarrow 1$ and $t\rightarrow \infty , $  gives $I=N$.
\subsection{\textbf{SIS Model}}
In SIS model, the population is divided into two disjoint classes that is susceptible to infectives and infectives to susceptibles.  The dynamics of the disease specified by two functions, the contact rate and the distribution of the infective period. Susceptible individuals become infective after contact with infective individuals. Infective individuals return to the susceptible class after an infective period. \\
In this model, we consider that the infected individual can recover and become susceptible again at a rate given by $\lambda I$, where $\lambda$ is positive constant. Thus, we get differential equations as  
\begin{align} \label{124}
{{dI} \over {dt}} &= \beta I\left( {N - I} \right) - \lambda I \\ {{dS} \over {dt}} &=  - \beta S\left( {N - S} \right) + \lambda I
\end{align}
\subsubsection{\textbf{Fractional Differential Equation  for SIS Model:}}
Here, we solve the given SIS model in the form of fractional differential equation.\\
From equation \eqref{124} we have,
\begin{align}\label{3.8}
{{dI} \over {dt}} =\beta I\left[ {A - I} \right], \hspace{1cm}  I\left( 0 \right) = {I_0},~A=\left( {N - {\lambda  \over \beta }} \right).
\end{align}
For developing fractional differential equation model we write \eqref{124} as,
\begin{align}\label{3.9}
\frac{d^\alpha I}{dt^\alpha} = \beta I\left( {A - I} \right), \hspace{1cm} I\left( 0 \right) = {I_0}.
\end{align}
As shown in the previous model, the solution of this fractional differential equation can also be given by 
\begin{align}\label{3.10}
{{{\log }_{{E_\alpha }}}I \ominus {{\log }_{{E_\alpha }}}\left( {A - I} \right)}&= {\frac{A\beta}{\Gamma(2-\alpha)}} \int t^{1-\alpha}dt^{\alpha} + c \\
\frac{I}{A-I}&=CE_\alpha \left[ {\frac{A\beta}{\Gamma(2-\alpha)}} \int t^{1-\alpha}dt^{\alpha} \right] ,
\end{align}
where $0 < \alpha  < 1$ and $\mathop {\lim }\limits_{\alpha  \to 1} \,\,{\log _{{E_\alpha }}}x = {\log _e}x$. 
Initially when time $t=0$, we have,
$I(0)=I_0$, and further simplification gives $C=\frac{I_0}{A-I_0}$.  
Thus, we arrive at,
\begin{align}\label{3.11}
& \frac{I}{A-I}=\frac{I_0}{A-I_0}E_\alpha \left[ {\frac{A\beta}{\Gamma(2-\alpha)}} \int t^{1-\alpha}dt^{\alpha} \right], \\ 
or \hspace{1cm}  &I=\frac{A}{1+\frac{A-I_0}{I_0} \left[ E_\alpha \lbrace {\frac{A\beta}{\Gamma(2-\alpha)}} \int t^{1-\alpha}dt^{\alpha} \rbrace \right]^{-1}}.
\end{align}
when $\alpha \rightarrow 1$ and $t\rightarrow \infty $ this yields $I=N-\frac{\lambda}{\beta},$ as reported by Kapur  \cite{kapur1988mathematical}.


\section{Conclusion}

Our approach for finding solution of fractional logistic equation is simple and may be useful in further studies in the field of mathematical modelling in fractional order and the theory of fractional calculus.
{\section*{Acknowledgement}
The first author is thankful to Sardar Vallabhbhai National Institute of Technology, Surat, Gujarat, INDIA for providing financial support in terms of Senior Research Fellowship.


\begin{landscape}

\begin{table}
\caption{Values of $\log_{E_\alpha}(x)$ for different values of $x$ and $\alpha $}.
\begin{tabular}{ | c | c | c | c | c | c | c | c | c | c | c | }
\hline
	$\alpha $ & 0.1 & 0.2 & 0.3 & 0.4 & 0.5 & 0.6 & 0.7 & 0.8 & 0.9 & 1 \\ \hline
$E_\alpha (=E_\alpha (1))$ & 23.1605 & 11.823 & 8.0407 & 6.1471 & 5.009 & 4.2486 & 3.7041 & 3.2946 & 2.9749 & 2.7183

\\ \hline
	$\log_{E_\alpha}(0.1)$ & -0.7327 & -0.9322 & -1.1046 & -1.2680 & -1.4291 & -1.5917 & -1.7584 & -1.9312 & -2.1120 & -2.3026 \\ \hline
	$\log_{E_\alpha}(0.2)$ & -0.5122 & -0.6516 & -0.7721 & -0.8863 & -0.9989 & -1.1126 & -1.2291 & -1.3499 & -1.4763 & -1.6094 \\ \hline
	$\log_{E_\alpha}(0.3)$ & -0.3831 & -0.4874 & -0.5776 & -0.6630 & -0.7472 & -0.8323 & -0.9194 & -1.0098 & -1.1043 & -1.2040 \\ \hline
	$\log_{E_\alpha}(0.4)$ & -0.2916 & -0.3710 & -0.4396 & -0.5046 & -0.5687 & -0.6334 & -0.6998 & -0.7685 & -0.8405 & -0.9163 \\ \hline
	$\log_{E_\alpha}(0.5)$ & -0.2206 & -0.2806 & -0.3325 & -0.3817 & -0.4302 & -0.4792 & -0.5293 & -0.5814 & -0.6358 & -0.6931 \\ \hline
	$\log_{E_\alpha}(0.6)$ & -0.1626 & -0.2068 & -0.2451 & -0.2813 & -0.3170 & -0.3531 & -0.3901 & -0.4284 & -0.4686 & -0.5108 \\ \hline
	$\log_{E_\alpha}(0.7)$ & -0.1135 & -0.1444 & -0.1711 & -0.1964 & -0.2214 & -0.2466 & -0.2724 & -0.2992 & -0.3272 & -0.3567 \\ \hline
	$\log_{E_\alpha}(0.8)$ & -0.0710 & -0.0903 & -0.1070 & -0.1229 & -0.1385 & -0.1542 & -0.1704 & -0.1872 & -0.2047 & -0.2231 \\ \hline
	$\log_{E_\alpha}(0.9)$ & -0.03353 & -0.0426 & -0.0505 & -0.05802 & -0.0654 & -0.0728 & -0.0805 & -0.0884 & -0.0966 & -0.1054 \\ \hline	

	$\log_{E_\alpha}(1)$ & 0 & 0 & 0 & 0 & 0 & 0 & 0 & 0 & 0 & 0 
\\ \hline
	$\log_{E_\alpha}(2)$ & 0.2206 & 0.2806& 0.3325 & 0.38169 & 0.4302 & 0.4792 & 0.5293 & 0.5814 & 0.6358 & 0.6931
\\ \hline
	$\log_{E_\alpha}(3)$ & 0.3496 & 0.4448 & 0.5270 & 0.6049 & 0.6818 & 0.7594 & 0.8390 & 0.9214 & 1.0077 & 1.0986
\\ \hline
	$\log_{E_\alpha}(4)$ & 0.4412 & 0.5612 & 0.6650 & 0.7634 & 0.8604 & 0.9583 & 1.0587 & 1.1627 & 1.2716 & 1.3863
\\ \hline
	$\log_{E_\alpha}(5)$ & 0.5122 & 0.6516 & 0.7721 & 0.8863 & 0.9989 & 1.1126 & 1.2291 & 1.3499 & 1.4763 & 1.6094
\\ \hline
	$\log_{E_\alpha}(6)$ & 0.5702 & 0.7254 & 0.8596 & 0.9867 & 1.1120 & 1.2386 & 1.3683 & 1.5028 & 1.6435 & 1.7917 
\\ \hline
	$\log_{E_\alpha}(7)$ & 0.6192 & 0.7878 & 0.9335 & 1.0715 & 1.2077 & 1.3452 & 1.4861 & 1.6321 & 1.7849 & 1.9459 
\\ \hline
	$\log_{E_\alpha}(8)$ & 0.6617 & 0.8419 & 0.9976 & 1.1451 & 1.2906 & 1.4375 & 1.5880 & 1.7441 & 1.9074 & 2.0794 
\\ \hline
	$\log_{E_\alpha}(9)$ & 0.6992 & 0.8895 & 1.0541 & 1.2099 & 1.3637 & 1.5189 & 1.6780 & 1.8429 & 2.0154 & 2.1972
\\ \hline
	$\log_{E_\alpha}(10)$ & 0.7327 & 0.9322 & 1.1046 & 1.2680 & 1.4291 & 1.5917 & 1.7584 & 1.9312 & 2.1120 & 2.3026 
\\ \hline

\end{tabular}
\label{table1}
 
\end{table}
\begin{table}
\caption{Examples of proposition 1 and 2.}
\begin{tabular}{|c|c|c|c|c|c|c|c|c|}

\hline 
$x_1$ & $x_2$ & $\alpha $ & $\log_{E_\alpha}(x_1 \cdot x_2)$ & $\log_{E_\alpha}\left(\frac{x_1}{ x_2}\right)$ & $\log_{E_\alpha}(x_1)$ & $\log_{E_\alpha}(x_2)$ & $\log_{E_\alpha}(x_1) +\log_{E_\alpha}(x_2)$ & $\log_{E_\alpha}(x_1) -\log_{E_\alpha}(x_2)$\\ 
\hline 
0.2 & 1 & 0.1 & -0.5122 & -0.5122 & -0.5122 & 0.00 & -0.5122 & -0.5122\\ 
\hline 
0.2 & 1 & 0.2 & -0.6516 & -0.6516 & -0.6516 & 0.00 & -0.6516 & -0.6516\\ 
\hline 
0.2 & 1 & 0.3 & -0.7721 & -0.7721 & -0.7721 & 0.00 & -0.7721 & -0.7721\\ 
\hline 
0.75 & 0.35 & 0.1 & -0.4256 & 0.2426 & -0.0915 & -0.3341 & -0.4256  & 0.2426 \\ 
\hline 
0.75 & 0.35 & 0.5 & -0.8301 & 0.4731 & -0.1785 & -0.6516 & -0.8301 & 0.4731\\ 
\hline 
0.75 & 0.35 & 0.9 & -1.2268 & 0.6991 & -0.2639 & -0.9630 & -1.2269 & 0.6991\\ 
\hline 
0.81 & 0.4 & 0.2 & -0.4563 & 0.2857 &  -0.0853 & -0.3709 & -0.4562 & 0.2856 \\ 
\hline 
0.81 & 0.4 & 0.7 & -0.8607 & 0.5388 & -0.1609 & -0.6998 & -0.8607 & 0.5389\\ 
\hline 
0.81 & 0.4 & 0.8 & -0.9453 & 0.5918 &-0.1767 & -0.7685 & -0.9453 & 0.5918 \\ 
\hline 
0.93 & 0.5 & 0.5 & -0.4752 & 0.3852 & -0.0450 & -0.4302 & -0.4752 & 0.3852 \\ 
\hline 
2 & 3 & 0.6 & 1.2386 & -0.2802& 0.4792 & 0.7594 & 1.2386 & -0.2802 \\ 
\hline 
3 & 5 & 0.7 & 2.0681 & -0.3901 & 0.8390 & 1.2291 & 2.0681 &  -0.3901\\ 
\hline 
6 & 7 & 0.8 & 3.1349 & -0.1293 & 1.5028 & 1.6321 & 3.1349 & -0.1293 \\ 
\hline 
10 & 2 & 0.9 & 2.7478 & 1.4763& 2.1121 & 0.6358 & 2.7479 & 1.4763\\ 
\hline 
\end{tabular} 
\label{table2}
\end{table}

\end{landscape}
\begin{figure}
\includegraphics[scale=.4]{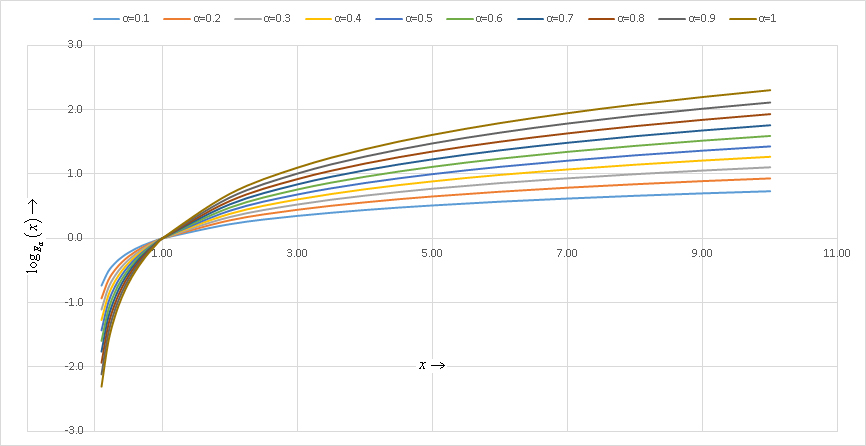}
\caption{Logarithmic function with base ${E_\alpha}$ or $E_\alpha(1)$ for different values of $x$ and $\alpha$.}
\label{Fig1}
\end{figure}

\end{document}